\documentclass[12pt]{amsart}
\usepackage{amsfonts}
\usepackage{ifthen}
\usepackage{amsthm}
\usepackage{amsmath}
\usepackage{graphicx}
\usepackage{amscd,amssymb,amsthm}

\newcounter{minutes}
\divide\time by 60
\newcounter{hours}
\multiply\time by 60 \addtocounter{minutes}{-\time}

\title[M. Arsenovi\'c, V. Manojlovi\'c, and  M. Vuorinen]
{H\"older continuity of harmonic quasiconformal mappings$^\bigstar$}

\author[]{Milo\v s Arsenovi\'c$\dagger$}
\author[]{Vesna Manojlovi\'c$\S$}
\author[]{Matti Vuorinen$\ddagger$}

\address{Department of Mathematics, University of Belgrade, Studentski trg 16, 11000 Belgrade, Serbia} \email{arsenovic@matf.bg.ac.rs}

\address{Faculty of Organizational Sciences, University of Belgrade, Jove Ilica 154, 11000 Belgrade, Serbia } \email{vesnam@fon.bg.ac.rs }

\address{Department of Mathematics, University of Turku, 20014 Turku,
Finland} \email{vuorinen@utu.fi}

\thanks{$\dagger$ Supported by Ministry of Science, Serbia, project M144010
\thanks{$\S$ Supported by Ministry of Science, Serbia, project M174024}
\thanks{$\ddagger$ Supported by the Academy of Finland, project 2600066611}
\thanks{$^\bigstar$ File:~\jobname .tex,
          printed: \number\year-\number\month-\number\day,
          \thehours.\ifnum\theminutes<10{0}\fi\theminutes} }

\newtheorem{theorem}[equation]{Theorem}

\newtheorem{lemma}[equation]{Lemma}

\newtheorem{corollary}[equation]{Corollary}

\numberwithin{equation}{section}

\begin{document}

\maketitle


\begin{abstract}
We prove that for harmonic quasiconformal mappings $\alpha$-H\"older continuity on the boundary implies $\alpha$-H\"older
continuity of the map itself. Our result holds for the class of uniformly perfect bounded domains, in fact we can allow
that a portion of the boundary is thin in the sense of capacity. The problem for general bounded domains remains open.
\end{abstract}

\section{Introduction}




The following theorem is the main result in \cite{MN}.

\begin{theorem}
Let $D$ be a bounded domain in $\mathbb R^n$ and let $f$ be a continuous mapping of $\overline D$ into $\mathbb R^n$
which is quasiconformal in $D$. Suppose that, for some $M > 0$ and $0 < \alpha \leq 1$,
\begin{equation} \label{Holder1}
|f(x) - f(y)| \leq M |x-y|^\alpha
\end{equation}
whenever $x$ and $y$ lie on $\partial D$. Then
\begin{equation}\label{Holder2}
|f(x) - f(y)| \leq M^\prime |x-y|^\beta
\end{equation}
for all $x$ and $y$ on $\overline D$, where
$\beta = \min (\alpha , K_I^{{1}/({1-n})} )$ and $M^\prime$ depends
only on  $M$, $\alpha$, $n$, $K(f)$ and $\mbox{\rm diam} (D)$.
\end{theorem}

The exponent $\beta$ is the best possible, as the example of a radial quasiconformal map $f(x) = |x|^{\alpha - 1} x$,
$0 < \alpha < 1$, of $\overline{\mathbb B^n}$ onto itself shows (see \cite{Va}, p. 49). Also, the assumption of boundedness is essential.
Indeed, one can consider $g(x) = |x|^a x$, $|x| \geq 1$ where $a > 0$. Then $g$ is quasiconformal in
$D = \mathbb R^n \setminus \overline{\mathbb B^n}$ (see \cite{Va}, p. 49), it is identity on $\partial D$ and hence
Lipschitz continuous on $\partial D$. However, $|g(te_1) - g(e_1)| \asymp t^{a+1}$, $t \to \infty$, and therefore $g$ is not globally
Lipschitz continuous on $D$.

This paper deals with the following question, suggested by P. Koskela: is it possible to replace $\beta$ with $\alpha$ if we assume, in addition to quasiconformality, that $f$ is harmonic? In the special case $D = \mathbb B^n$ this was proved, for arbitrary moduli of continuity $\omega(\delta)$, in \cite{ABM}. Our main result is that the answer is positive, if $\partial D$ is a uniformly perfect set (cf. \cite{JV}). In fact, we prove a more general result, including domains having a thin, in the sense of capacity, portion of the boundary. However, this generality is in a sense illusory, because any harmonic and quasiconformal (briefly hqc) mapping extends harmonically and quasiconformally across such portion of the boundary. Nevertheless, it leads to a natural open question: is the answer positive for arbitrary bounded domain in $\mathbb R^n$?

In the case of smooth boundaries much better regularity up to the boundary can be deduced, see \cite{Ka};
related results for harmonic functions were obtained by \cite{Ai}.

We denote by $B(x, r)$ and $S(x, r)$ the open ball, respectively sphere, in $\mathbb R^n$ with center $x$ and radius $r>0$.
We adopt the basic notation, terminology and definitions related to quasiconformal maps from \cite{Va}. A condenser is a pair $(K, U)$,
where $K$ is a non-empty compact subset of an open set $U \subset \mathbb R^n$. The capacity of the condenser $(K, U)$ is defined as
$$ \mbox{\rm cap} (K, U) = \inf \int_{\mathbb R^n} | \nabla u |^n dV,$$
where infimum is taken over all continuous real-valued $u \in ACL^n (\mathbb R^n)$ such that $u(x) = 1$ for $x \in K$ and $u(x) = 0$ for
$x \in \mathbb R^n \setminus U$. In fact, one can replace the $ACL^n$ condition with Lipschitz continuity in this definition.  We note that,
for a compact $K \subset \mathbb R^n$ and open bounded sets $U_1$ and $U_2$ containing $K$ we have: $\mbox{\rm cap}(K, U_1) = 0$ iff
$\mbox{\rm cap}(K, U_2) = 0$, therefore the notion of a compact set of zero capacity is well defined (see \cite{Vu}, Remarks 7.13) and we
can write $\mbox{\rm cap}(K) = 0$ in this situation. For the notion of the modulus $M(\Gamma)$ of a family $\Gamma$ of curves in
$\mathbb R^n$ we refer to \cite{Va} and \cite{Vu}. These two notions are related: by results of \cite{He} and \cite{Zi} we have
$$ \mbox{\rm cap}(K, U) = M(\Delta(K, \partial U; U)),$$
where $\Delta(E, F; G)$ denotes the family of curves connecting $E$ to $F$ within $G$, see \cite{Va} or \cite{Vu} for details.

In addition to this notion of capacity, related to quasiconformal mappings, we need Wiener capacity, related to harmonic functions. For a compact $K \subset \mathbb R^n$, $n \geq 3$, it is defined by
$$\mbox{\rm cap}_W (K) = \inf \int_{\mathbb R^n} |\nabla u|^2 dV,$$
where infimum is taken over all Lipschitz continuous compactly supported functions $u$ on $\mathbb R^n$ such that $u = 1$ on $K$. Let us note that every compact $K \subset \mathbb R^n$ which has capacity zero has Wiener capacity zero. Indeed, choose an open ball $B_R = B(0, R) \supset K$. Since $n \geq 2$ we have, by H\"older inequality,
$$\int_{\mathbb R^n} |\nabla u|^2 dV \leq |B_R|^{1 - 2/n} \left( \int_{\mathbb R^n} |\nabla u|^n dV \right)^{2/n}$$
for any Lipschitz continuous $u$ vanishing outside $U$, our claim follows immediately from definitions.

A compact set $K \subset \mathbb R^n$, consisting of at least two points, is $\alpha$-uniformly perfect ($\alpha > 0$) if there is no ring $R$
separating $K$ (i.e. such that both components of $\mathbb R^n \setminus R$ intersect $K$) such that $\mbox{\rm mod}(R) > \alpha$. We say
that a compact $K \subset \mathbb R^n$ is uniformly perfect if it is $\alpha$-uniformly perfect for some $\alpha > 0$.

We denote the $\alpha$-dimensional Hausdorff measure of a set $F \subset \mathbb R^n$ by $\Lambda_\alpha (F)$.





\section{The main result}

In this section $D$ denotes a bounded domain in $\mathbb R^n$, $n \geq 3$. Let
$$\Gamma_0 = \{ x \in \partial D : \mbox{\rm cap}\, (\overline B(x, \epsilon) \cap \partial D) = 0 \; \mbox{\rm for some}\; \epsilon > 0 \},$$
and $\Gamma_1 = \partial D \setminus \Gamma_0$. Using this notation we can state our main result.

\begin{theorem} \label{my21}
Assume $f : \overline D \rightarrow \mathbb R^n$ is continuous on $\overline D$, harmonic and quasiconformal in $D$. Assume
$f$ is H\"older continuous with exponent $\alpha$, $0 < \alpha \leq 1$, on $\partial D$ and $\Gamma_1$ is uniformly perfect.
Then $f$ is H\"older continuous with exponent $\alpha$ on $\overline D$.
\end{theorem}

If $\Gamma_0$ is empty we obtain the following

\begin{corollary} \label{myc22}
If $f : \overline D \rightarrow \mathbb R^n$ is continuous on $\overline D$, H\"older continuous with exponent $\alpha$,
$0 < \alpha \leq 1$, on $\partial D$, harmonic and quasiconformal in $D$ and if $\partial D$ is uniformly perfect, then
$f$ is H\"older continuous with exponent $\alpha$ on $\overline D$.
\end{corollary}

The first step in proving Theorem \ref{my21} is reduction to the case $\Gamma_0 = \emptyset$. In fact, we show that existence of a hqc
extension of $f$ across $\Gamma_0$ follows from well known results. Let $D^\prime = D \cup \Gamma_0$. Then
$D^\prime$ is an open set in $\mathbb R^n$, $\Gamma_0$ is a closed subset of $D^\prime$ and $\partial D^\prime = \Gamma_1$.

Clearly $\mbox{\rm cap} (K \cap \Gamma_0) = 0$ for each compact $K \subset D^\prime$, and therefore, by Lemma 7.14 in \cite{Vu}, $\Lambda_\alpha (K \cap \Gamma_0) = 0$ for each $\alpha > 0$. In particular, $\Gamma_0$ has $\sigma$-finite $(n-1)$-dimensional Hausdorff measure. Since it is closed in $D^\prime$, we can apply Theorem 35.1 in \cite{Va} to conclude that $f$ has a quasiconformal extension $F$ across $\Gamma_0$ which has the same quasiconformality constant as $f$.

Since $\Gamma_0$ is a countable union of compact subsets $K_j$ of capacity zero and hence of Wiener capacity zero we
conclude that $\Gamma_0$ has Wiener capacity zero. Hence, by a classical result (see \cite{Car}), there is a (unique) extension $G : \overline{D^\prime} \rightarrow \mathbb R^n$ of $f$ which is harmonic in $D^\prime$. Obviously, $F = G$ is a harmonic quasiconformal extension of $f$ to $\overline{D^\prime}$ which has the same quasiconformality constant as $f$.

In effect, we reduced the proof of Theorem \ref{my21} to the proof of Corollary \ref{myc22}. We begin the proof of Corollary \ref{myc22}
with the following

\begin{lemma} \label{myl23}
Let $D \subset {\mathbb R}^n$ be a bounded domain with uniformly perfect boundary.
There exists a constant $m > 0$ such that for every $y \in D$ we have
\begin{equation}\label{cap}
\mbox{\rm cap}( \overline B(y, \frac{d}{2}), D) \geq m \;\; ,\;\;\;\;\;\;\;\;\;\; d = \mbox{\rm dist}(y, \partial D).
\end{equation}
\end{lemma}

{\it Proof.} 
Fix $y \in D$ as above and $z \in \partial D$ such that
$|y-z| = d \equiv r \,.$ Clearly $diam(\partial D) = diam(D)> 2r\,.$
Set $F_1 = \overline{B}(z,r) \cap (\partial D)$
and $F_2 = \overline{B}(z,r) \cap \overline{B}(y, \frac{d}{2})$, $F_3 = S(z,2r) \,.$ Let $\Gamma_{i,j} = \Delta(F_i,F_j; \mathbb{R} ^n)$ for
$i,j=1,2,3$. By \cite[Thm 4.1(3)]{JV} there exists a constant $a=a(E,n)>0$ such that
$$
M(\Gamma_{1,3}) \ge a
$$
while by standard estimates \cite[7.5]{Va} there exists $b=b(n)>0$ such that
$$
M(\Gamma_{2,3}) \ge b \,.
$$
Next, by \cite[Cor 5.41]{Vu} there exists $m=m(E,n)>0$ such that
$$
M(\Gamma_{1,2}) \ge m \,.
$$
Finally, with $B= {\overline B}(y, d/2)$ we have
$$ {\rm cap} (B, D)= M(\Delta(B, \partial D;  \mathbb{R} ^n)) \ge M(\Gamma_{1,2}) \ge m\,. $$
\hfill $\square$

\bigskip

In conclusion,
from the above lemma, our assumption
$$|f(x_1) - f(x_2)| \leq C|x_1 - x_2|^\alpha \;, \;\;\;\;\;\;\;\; x_1, x_2 \in \partial D,$$
and Lemma 8 in \cite{MN} we conclude that there is a constant $M$, depending on $m$, $n$, $K(f)$, $C$ and $\alpha$ only such that
$$|f(x) - f(y)| \leq M |x - y|^\alpha \;, \;\;\; y \in D,\; x \in \partial D,\; \mbox{\rm dist}(y, \partial D) = |x-y|.$$
However, an argument presented in \cite{MN} shows that the above estimate holds for $y \in D$, $x \in \partial D$ without any
further conditions, but with possibly different constant:
\begin{equation}\label{bdryest}
|f(x) - f(y)| \leq M^\prime |x - y|^\alpha \;, \;\;\; y \in D,\; x \in \partial D.
\end{equation}

The following lemma was proved in \cite{CK} for
real valued functions, but the proof relies on the maximum principle which holds also for vector valued harmonic functions,
hence lemma holds for harmonic mappings as well.

\begin{lemma} \label{myl24}
Assume $h : \overline D \rightarrow \mathbb R^n$ is continuous on $\overline D$ and harmonic in $D$. Assume
for each $x_0 \in \partial D$ we have
$$ \sup_{B_r(x_0) \cap D^\prime} | h(x) - h(x_0) | \leq \omega(r) \;\;\;\;\;\;\;\;\; \mbox{ for } \;\; 0 < r \leq r_0.$$
Then $|h(x) - h(y)| \leq \omega(|x-y|)$ whenever $x, y \in D$ and $|x-y| \leq r_0$.
\end{lemma}

Now we combine (\ref{bdryest}) and the above lemma, with $r_0 = \mbox{\rm diam}(D)$, to complete the proof of Corollary \ref{myc22} and therefore
of Theorem \ref{my21} as well.

\end{document}